\documentclass[12pt,a4paper]{article}
\usepackage[dvips]{graphicx}
\usepackage{latexsym,amsfonts,amscd,subfigure}
\usepackage{amsmath,amssymb,amsthm,mathrsfs,dsfont}
\usepackage{epsfig}
\usepackage{changebar}
\usepackage{indentfirst}
\usepackage{verbatim}
\usepackage{color}
\usepackage{amsmath}
\usepackage{amssymb}

\usepackage{multirow}
\usepackage{epstopdf}
\graphicspath{{figures/}}

\numberwithin{equation}{section}

\begin{document}

\baselineskip=2pc

\newtheorem{algorithm}{ALGORITHM}
\newtheorem{remark}{REMARK}
\newtheorem{theorem}{THEOREM}
\newtheorem{proposition}{PROPOSITION}
\newtheorem{example}{EXAMPLE}

\title{Finite volume HWENO schemes for nonconvex conservation laws\footnote{Research
was supported by NSFC grants 11571290, 91530107, Air Force Office of Scientific Research
FA9550-16-1-0179 and NSF DMS-1522777.}}

\author{Xiaofeng Cai\footnote{Department of Mathematics, University of Houston, Houston, 77204.
E-mail: xfcai@math.uh.edu.},
Jianxian Qiu\footnote{School of
Mathematical Sciences and Fujian Provincial Key Laboratory of Mathematical
Modeling \&
High-Performance Scientific Computing, Xiamen University, Xiamen, Fujian,
361005, P.R. China.
E-mail: jxqiu@xmu.edu.cn.},
and Jingmei Qiu\footnote{Department of Mathematics, University of Houston, Houston, 77204. E-mail: jingqiu@math.uh.edu.}}

\date{}
\maketitle
{Abstract:
Following the previous work of Qiu and Shu [SIAM J. Sci. Comput., 31 (2008), 584-607], we investigate the performance of Hermite weighted essentially non-oscillatory (HWENO) scheme for nonconvex conservation laws. Similar to many other high order methods, we show that the finite volume HWENO scheme performs poorly for some nonconvex conservation laws. We modify the scheme around the nonconvex regions, based on a first order monotone scheme and a second entropic projection, to ensure entropic convergence. Extensive numerical tests are performed. Compare with the earlier work of Qiu and Shu which focuses on 1D scalar problems, we apply the modified schemes (both WENO and HWENO) to one-dimensional Euler system with nonconvex equation of state and two-dimensional problems.
}


\newpage
\section{Introduction}

In this paper, we consider the Cauchy problem for nonconvex hyperbolic  conservation laws:
\begin{equation}
\label{EQ}
\left\{
\begin{array}{c}
q_t+\nabla \cdot F(q)=0, \\
q(\mathbf{x}, 0)=q_0(\mathbf{x}),
\end{array}
\right.
\end{equation}
whose entropy solutions may admit composite waves
which involve a sequence of shocks and rarefaction waves and are difficult to be resolved numerically.
Such examples include scalar conservation laws with nonconvex flux functions and hyperbolic systems such as the Euler system and magnetohydrodynamics system with a nonconvex equation of state (EOS) \cite{menikoff1989riemann,wang1999second,muller2006riemann,serna2014anomalous}.

It is well known that first order monotone schemes converge to entropy solutions of both convex and nonconvex conservation laws \cite{crandall1980monotone}, but with a relatively slow convergence rate. It has also been known \cite{kurganov2007adaptive, qiu_siam_sc} that there are some nonconvex conservation laws, for which high order schemes such as the ones with weighted essentially non-oscillatory (WENO) reconstructions \cite{weno_review} and discontinuous Galerkin methods \cite{cockburn1989tvb} would fail to converge to the entropy solution.
There have been great research effort in ensuring entropic convergence for general nonlinear conservation laws, for example by adding entropy viscosity \cite{guermond2011entropy} and by modifying reconstruction operators. Examples for the latter approach include the computationally inexpensive strategy proposed in \cite{kurganov2007adaptive} on an adaptive choice between a low order dissipative reconstruction and a high order central WENO scheme, as well as low order modifications around nonconvex regions to ensure entropic convergence proposed in \cite{qiu_siam_sc}.  Compare with the work in \cite{kurganov2007adaptive}, with more computational effort, the second order entropic convergence of the schemes can be rigorously proved \cite{bouchut1996muscl}.
%

This paper is a natural extension of our earlier work in \cite{qiu_siam_sc}. We investigate the performance of the finite volume Hermite WENO (HWENO) scheme for nonconvex conservation laws and apply the corresponding modifications as being done in \cite{qiu_siam_sc}. In addition to the scalar examples discussed in \cite{qiu_siam_sc}, we investigate the performance of modified WENO and HWENO scheme for  2D problems, as suggested in \cite{kurganov2007adaptive}. The FV HWENO scheme was originally proposed in \cite{qiu2004hermite, qiu2005hermite}.
The key idea of the scheme is to evolve more pieces of information, i.e. functions and their spatial gradients, per computational cell. With such mechanism, the HWENO scheme has relatively compact stencils, hence it is easier to handle boundary conditions compared with the traditional WENO scheme \cite{weno_review}. Moreover, with the same formal accuracy, compact stencils are known to exhibit better resolution of small scale structures by improving dispersive and dissipative properties.

%

An outline of this paper is as follows.
Section \ref{HermiteWENO} describes the high order FV HWENO scheme.
In Section \ref{bad_example}, FV HWENO schemes with a first order monotone modification and a second order modification using an entropic projection around nonconvex regions are proposed for nonconvex conservation laws.
In Section \ref{numericaltests}, numerical examples are shown to demonstrate the effectively of proposed schemes.
Concluding remarks are given in Section \ref{concluding}.

\section{Description of FV HWENO schemes}
\label{HermiteWENO}
We briefly review the  FV HWENO scheme for
solving conservation laws
 \cite{qiu2004hermite, qiu2005hermite, zhu2008class}.
The idea of HWENO method is to numerically evolve {\em both the function and its spatial gradients}, and use these information in the reconstruction process. Thus it leads to a more compact reconstruction stencil compared with the traditional WENO scheme \cite{shi2002technique,weno_review}.

\noindent
{\underline{\bf General scheme formulation of HWENO.}}
Taking the gradient with respect to spatial variables in \eqref{EQ}, we obtain the evolution equation for function's gradients,
\begin{equation}
(\nabla q)^T_t +\nabla \cdot( \nabla  \otimes F(q) ) =0,
\end{equation}
where $\otimes$ is a tensor product.
The FV HWENO scheme is defined for the equations:
\begin{equation}
\mathbf{U}_t + \nabla \cdot \mathcal{F} (\mathbf{U})=0,
\label{hermite_eq}
\end{equation}
where $\mathbf{U}=(q,\nabla q)^T$ and $\mathcal{F}(\mathbf{U})=\left(
\begin{array}{c}
F(q) \\
\nabla  \otimes F(q)
\end{array}\right)$.
We integrate the system \eqref{hermite_eq} on a control volume $\Omega_k$, which
is an interval $I_j=[x_{j-\frac12},x_{j+\frac12}]$ for 1D cases or a
rectangle $[x_{i-\frac12},x_{i+\frac12}]\times[y_{j-\frac12},y_{j+\frac12}]$ for
2D cases.
The integral form of \eqref{hermite_eq} reads,
\begin{equation}
\frac{d}{dt}\overline{\mathbf{U}}_{\Omega_k} = -\frac{1}{|\Omega_k|} \int_{\partial
\Omega_k} \mathcal{F} (\mathbf{U})\cdot \mathbf{n} ds
\label{semi_hermite}
\end{equation}
where $|\Omega_k|$ is the volume of  $\Omega_k$ and $\mathbf{n}$ represents the outward unit normal vector to  $\partial
\Omega_k$.
The line integral in \eqref{semi_hermite} can be approximated  by a
$L$-point Gaussian quadrature on each side of ${\partial
\Omega_k}=\bigcup_{s=1}^S \partial \Omega_{ks}$:
\begin{equation}
\int_{\partial \Omega_k} \mathcal{F}(\mathbf{U})\cdot \mathbf{n} ds \approx  \mathop{\sum}_{s=1}^S
|\partial \Omega_{ks}| \mathop{\sum}_{l=1}^L \omega_l \mathcal{F}(\mathbf{U}(G_{sl}, t) )\cdot \mathbf{n},
\label{semi_hermite1}
\end{equation}
where $G_{{sl}}$ and $\omega_l$ are Gaussian quadrature points on $\partial \Omega_{ks}$ and weights
respectively. $\mathcal{F}(\mathbf{U}(G_{sl}, t) )\cdot \mathbf{n}$
is evaluated by a numerical flux (approximate or exact Riemann
solvers).
We adopt the Lax-Friedrichs flux in this paper, which is given by
\begin{equation*}
\mathcal{F}(\mathbf{U}(G_{sl}, t))\cdot \mathbf{n} \approx \frac12[ \mathcal{F}(\mathbf{U}^-(G_{sl}, t) ) +
\mathcal{F}(\mathbf{U}^+(G_{sl}, t) ) ] \cdot \mathbf{n} -\alpha (\mathbf{U}^+(G_{sl}, t) - \mathbf{U}^-(G_{sl}, t)
),
\label{semi_hermite2}
\end{equation*}
where $\alpha$ is taken as an upper bound for eigenvalues of the Jacobian along the direction $\mathbf{n}$,
and $\mathbf{U}^-$ and $\mathbf{U}^+$   are the reconstructed values of $\mathbf{U}$ at Gaussian point $G_{sl}$ inside and outside $\Omega_{k}$.
Finally, the semi-discretization
HWENO scheme (\ref{semi_hermite}) can be written in the following ODE form:
\begin{equation}
\frac{d}{dt}\overline{\mathbf{U}} = \mathcal{L}(\overline{\mathbf{U}}).
\label{fv_ODE}
\end{equation}
The  ODE system \eqref{fv_ODE} is then discretized in time by a strong stability preserving
Runge-Kutta (RK) method in \cite{shu1988efficient}. The following third-order version is
used in this paper,
\begin{equation}
\begin{array}{l}
\overline{\mathbf{U}}^{(1)}=\overline{\mathbf{U}}^n+\Delta t \mathcal{L}(\overline{\mathbf{U}}^n) , \\
\overline{\mathbf{U}}^{(2)}=\frac{3}{4}\overline{\mathbf{U}}^n+\frac{1}{4}(\overline{\mathbf{U}}^{(1)}
+\Delta t \mathcal{L}(\overline{\mathbf{U}}^{(1)})),\\
\overline{\mathbf{U}}^{n+1}=\frac{1}{3}\overline{\mathbf{U}}^n+\frac{2}{3}(\overline{\mathbf{U}}^{(2)}
+\Delta t \mathcal{L}(\overline{\mathbf{U}}^{(2)})).
\end{array}
\label{SSP3}
\end{equation}

\noindent
{\underline{\bf A scalar 1D example.}} As an example, we consider a scalar 1D equation,
 \begin{equation}
 \label{eq:1dconservationlaws}
 q_t + f(q)_x =0.
 \end{equation}
Taking the derivative of \eqref{eq:1dconservationlaws}, we obtain the equation for the derivative,
 \begin{equation}
 \label{derivative}
 \xi_t + \mathcal{H}(q,\xi)_x = 0,
 \end{equation}
where $\xi = q_x$ and $\mathcal{H}(q,\xi) = f'(q)\xi.$
Let $\overline{q}_j$ and $\overline{\xi}_j$ denote approximation to cell averages of $q$ and  $\xi$ over cell $I_j$ respectively, the semi-discrete FV HWENO scheme is designed by approximating spatial derivatives in equation \eqref{eq:1dconservationlaws} and \eqref{derivative} with the following flux difference form,
\begin{equation}
\label{label:scheme}
\begin{cases}
\frac{d \overline{q}_j }{ dt } = -\frac{  \widehat{f}\left(q_{j+\frac12}^-,q_{j+\frac12}^+\right) - \widehat{f}\left(q_{j-\frac12}^-,q_{j-\frac12}^+\right)  }{\Delta x}   ,\\
\frac{d \overline{\xi}_j }{dt} = -\frac{ \widehat{\mathcal{H}} \left( q_{j+\frac12}^-,q_{j+\frac12}^+;\xi_{j+\frac12}^-,\xi_{j+\frac12}^+ \right) -
\widehat{\mathcal{H}} \left( q_{j-\frac12}^-,q_{j-\frac12}^+;\xi_{j-\frac12}^-,\xi_{j-\frac12}^+ \right)
 }{\Delta x},
\end{cases}
\end{equation}
where $q_{j+\frac12}^\pm$ and $\xi_{j+\frac12}^\pm$ are reconstructed with high order from neighboring cell averages $\overline{q}$ and $\overline{\xi}$. The details of such reconstruction procedures can be found in \cite{qiu2004hermite}. $\widehat{f}(a,b)$ is a monotone numerical flux (non-decreasing in the first argument and non-increasing in the second argument), and $\widehat{\mathcal{H}}(a,b;c,d)$ is non-decreasing in the third argument and non-increasing in the fourth argument.
In this paper, we use the Lax-Friedrichs flux in \cite{qiu2004hermite},
\begin{equation}
\label{llf}
\begin{split}
\widehat{f}(a,b)=\frac12[f(a)+f(b)-\alpha(b-a)],\\
\widehat{\mathcal{H}}(a,b;c,d)=\frac12[\mathcal{H}(a,c)+\mathcal{H}(b,d)-\alpha(d-c)],
\end{split}
\end{equation}
where $\alpha=\mathop{\max}_{q}|f'(q)|$. For the first order ``building block" of the HWENO scheme with the Lax-Friedrichs flux, the total variation stability is proved in \cite{qiu2004hermite}.

\noindent
{\underline{\bf A scalar 2D example.}}
We consider a 2D problem on a rectangular domain $[a,b]\times[c,d]$:
\begin{equation}
q_t + f(q)_x + g(q)_y = 0.
\label{2dCL}
\end{equation}
We consider a set of uniform mesh with $I_{ij} = [x_{i-\frac12},x_{i+\frac12} ]\times[y_{j-\frac12},y_{j+\frac12} ]$,
\begin{equation*}
a = x_{\frac12} < x_{\frac32} < \cdots < x_{N_x-\frac12} <x_{N_x+\frac12}=b,\ \Delta x =\frac{b-a}{N_x},
\end{equation*}
\begin{equation*}
c= y_{\frac12} < y_{\frac32} < \cdots < y_{N_y - \frac12} < y_{N_y+\frac12} = d,\  \Delta y =\frac{d-c}{N_y}.
\end{equation*}
Let $\xi=\frac{\partial q}{\partial x}$, $\eta=\frac{\partial q}{\partial y}$ and $\overline{q}_{ij}=\frac{1}{\Delta x \Delta y} \int_{I_{ij} } q dx dy$,  $\overline{\xi}_{ij}=\frac{1}{\Delta y} \int_{I_{ij} } \frac{\partial q}{\partial x} dx dy$, $\overline{\eta}_{ij}=\frac{1}{\Delta x} \int_{I_{ij} } \frac{\partial q}{\partial y} dx dy$ be cell averages.
Taking spatial derivatives of \eqref{2dCL}, we obtain
\begin{align}
\label{2d1}
\xi_t + \mathcal{H}(q,\xi)_x + \mathcal{R}(q,\xi)_y =0,\\
\eta_t + \mathcal{K}(q,\eta)_x + \mathcal{S}(q,\eta)_y =0,
\label{2d2}
\end{align}
where $\mathcal{H}(q,\xi)=f'(q)\xi, \mathcal{R}(q,\xi)= g'(q)\xi,
\mathcal{K}(q,\eta)=f'(q)\eta,\mathcal{S}(q,\eta)=g'(q)\eta$.
A semi-discrete FV HWENO discretization is given by
\begin{equation}
\begin{cases}
\frac{d}{dt}\overline{q}_{ij} = -\frac{1}{\Delta x}(\widehat{f}_{i+\frac12,j} - \widehat{f}_{i-\frac12,j} )
- \frac{1}{\Delta y}(\widehat{g}_{i,j+\frac12} -\widehat{g}_{i,j-\frac12} ),
\\[3mm]
\frac{d}{dt}\overline{\xi}_{ij} = -\frac{1}{\Delta x}(\widehat{ \mathcal{H} }_{i+\frac12,j} - \widehat{ \mathcal{H} }_{i-\frac12,j} )
- \frac{1}{\Delta y}(\widehat{  \mathcal{R} }_{i,j+\frac12} -\widehat{  \mathcal{R} }_{i,j-\frac12} ),
\\[3mm]
\frac{d}{dt}\overline{\eta}_{ij} = -\frac{1}{\Delta x}(\widehat{ \mathcal{K} }_{i+\frac12,j} - \widehat{ \mathcal{K} }_{i-\frac12,j} )
- \frac{1}{\Delta y}(\widehat{  \mathcal{S} }_{i,j+\frac12} -\widehat{  \mathcal{S} }_{i,j-\frac12} ).
\end{cases}
\label{2dsemi}
\end{equation}
We define
\begin{align}
\label{eq: flux_r}
\widehat{f}_{i+\frac12,j} = \frac{1}{\Delta y} \int_{y_{j-\frac12} }^{y_{j+\frac12} } f(q( x_{i+\frac12} ,y) ) dy \approx \mathop{\sum_{i_g}^L } \omega_{i_g} \widehat{f}(q_{i+\frac12,i_g}^-,q_{i+\frac12,i_g}^+ ),\\
\widehat{\mathcal{H}}_{i+\frac12,j} = \frac{1}{\Delta y} \int_{y_{j-\frac12} }^{y_{j+\frac12} } f'(q( x_{i+\frac12} ,y) )\xi dy
\approx \mathop{\sum_{i_g}^L } \omega_{i_g} \widehat{H}(q_{i+\frac12,i_g}^-,q_{i+\frac12,i_g}^+, \xi_{i+\frac12,i_g}^-,\xi_{i+\frac12,i_g}^+), \\
\widehat{ \mathcal{K} }_{i+\frac12,j} = \frac{1}{\Delta y} \int_{y_{j-\frac12} }^{y_{j+\frac12} } f'(q( x_{i+\frac12} ,y) )\eta dy
\approx \mathop{\sum_{i_g}^L } \omega_{i_g} \widehat{K}(q_{i+\frac12,i_g}^-,q_{i+\frac12,i_g}^+, \eta_{i+\frac12,i_g}^-,\eta_{i+\frac12,i_g}^+),
\end{align}
as the average of fluxes over the right boundary of cell $I_{ij}$, where the integrations are evaluated by applying the $L$-point Gauss quadrature.
The flux functions $\widehat{f}$, $\widehat{H}$ and $ \widehat{K}$ are taken as the Lax-Friedrich flux as in the 1D, see eq.~\eqref{llf}, and $q_{i+\frac12,i_g}^\pm$ are reconstructed with high order by neighboring cell averages of $q$, $\xi$ and $\eta$. In the next paragraph, we briefly describe such reconstruction procedure. $\widehat{g}_{i,j+\frac12}$, $\widehat{\mathcal{R}}_{i,j+\frac12}$ and $\widehat{\mathcal{S}}_{i,j+\frac12}$ in eq.~\eqref{2dsemi} are evaluated in a similar fashion as the average of fluxes over the top boundary of a cell.

\begin{figure}[h!]
\centering
\includegraphics[height=70mm]{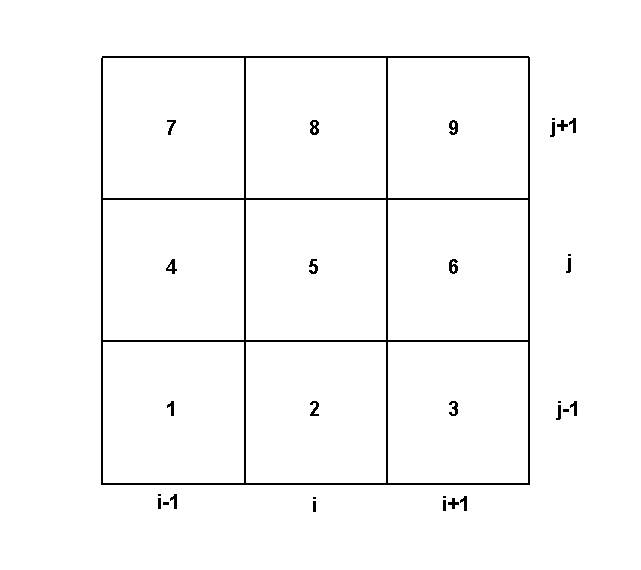}
\caption{\small The big stencil.}
\label{big}
\end{figure}

We only review the fourth order reconstruction in 2D and refer to \cite{zhu2008class} for more details. We relabel the cell $I_{ij}$ and its neighboring cells as $I_1,\cdots,I_9$ as shown in Figure \ref{big}, where $I_{ij}$ is relabeled as $I_5$. We construct the quadratic polynomials $p_n(x,y)$ $(n=1,\cdots,8)$ in the following stencils, $S_1 = \{ I_1 , I_2 ,I_4, I_5 \}$,
$S_2 = \{I_2,I_3,I_5,I_6\},$
$S_3 = \{ I_4, I_5, I_7, I_8 \},$
$S_4 = \{ I_5, I_6, I_8, I_9 \},$
$S_5 = \{ I_1,I_2,I_3, I_4,I_5,I_7 \},$
$S_6 = \{ I_1,I_2,I_3, I_5,I_6,I_9 \},$
$S_7 = \{ I_1,I_4,I_5, I_7,I_8,I_9 \},$
$S_8 = \{ I_3,I_5,I_6, I_7,I_8,I_9 \}$ to approximate $q(x,y)$. For instance, a quadratic polynomial can be reconstructed based on the information $\{\overline{q}_1,\overline{q}_2,\overline{q}_4,\overline{q}_5,\overline{\xi}_4,\overline{\eta}_2\}$ in the stencil $S_1$.
Such reconstruction will reconstruct a quadratic polynomial on $I_{ij}$. Similar reconstructions can be done for stencil $S_2$, $S_3$ and $S_4$.
For stencil $S_5$ to $S_8$, only cell averages are used in the reconstruction process. We remark that other combination of information are possible for reconstructing 2D quadratic polynomial. The one we just mentioned seems to be very robust and is  implemented in our numerical experiments.
 If we choose the linear weights denoted by $\gamma_{1}^{(l)},\cdots,\gamma_8^{(l)}$ such that
\begin{equation}
q(G_l) = \mathop{\sum_{n=1}^{8} }\gamma_n^{(l)} p_n(G_l)
\label{4th}
\end{equation}
is valid for any polynomial $q$ of degree at most 3, leading to  a fourth-order approximation of $q$ at the point $G_l$ for all sufficiently smooth functions $q$.
Notice that \eqref{4th} holds for any polynomial $q$ of degree at most 2 if $ \mathop{\sum_{n=1}^8}\gamma_n^{(l)}=1$.
There are four additional  constraints on the linear weights $\gamma_1^{(l)},\cdots,\gamma_8^{(l)}$ so that  \eqref{4th}  holds for $q=x^3,x^2y,xy^2$ and $y^3$. The rest of free parameters are determined by  a least square procedure to minimize $\mathop{\sum_{n=1}^8 } (\gamma_n^{(l)})^2 $.

As for the derivatives (e.g.  $\xi^-(G_l,t)$), a third-order approximation in each stencil is enough to obtain the fourth-order approximation to $q(x,y)$. For instance, a  cubic polynomial on $I_5$ can be reconstructed based on the information $\{\overline{q}_1,$$\overline{q}_2,$$\overline{q}_4,$$\overline{q}_5,$$\overline{\xi}_1,$$\overline{\xi}_4,$  $\overline{\xi}_5,$$\overline{\eta}_1, $$\overline{\eta}_2, $$\overline{\eta}_5  \}$ in the stencil $S_1$. Similar reconstructions can be done for stencil $S_2$, $S_3$ and $S_4$. The information $\{\overline{\xi}_1,\overline{\xi}_2,\overline{\xi}_3,\overline{\xi}_4,\overline{\xi}_5,\overline{\xi}_7\}$ in the stencil $S_5$ is adopted to approximate $\xi^-(G_l,t)$. Finally, $\gamma_l=\frac18\ (l=1,\cdots,8)$ can be chosen.
 The  nonlinear weights of 2D HWENO reconstruction can be designed by following the way of the WENO method.

\section{The modified FV HWENO schemes for nonconvex conservation laws}
\label{bad_example}

Although FV HWENO schemes can be successfully applied in many applications \cite{qiu2004hermite,qiu2005hermite,zhu2008class,HWENOHJ,zheng2016directly},  they perform poorly for some nonconvex conservation laws as shown below. To remedy this, we propose a first order monotone modification and a second order modification with an entropic projection around nonconvex regions.

\subsection{An example of nonconvex conservation laws with poor performance for the FV HWENO scheme}
\label{ssub: bad_example}
\begin{figure}[h!]
\centering
\includegraphics[height=60mm]{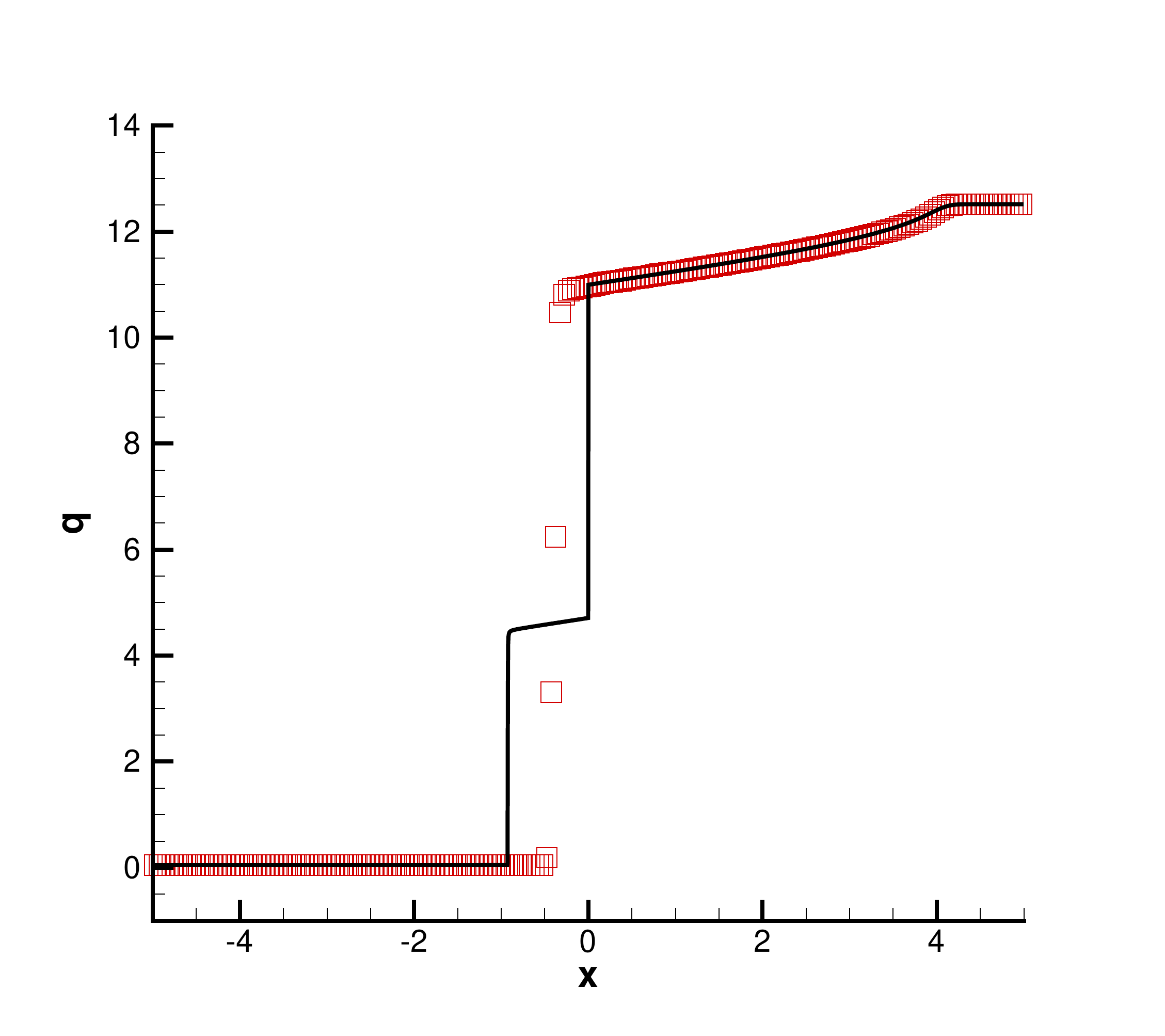}
\caption{\small Solid lines: the reference solution of \eqref{eq:1dconservationlaws} at the time $t=4$; Squares: FV HWENO scheme with Lax-Friedrichs flux and a uniform mesh $\Delta x = 0.05$. CFL=0.5.}
\label{Fig:sin1}
\end{figure}

We first show a nonconvex conservation law, for which the FV HWENO scheme performs poorly  in converging to the entropy solution. We consider the scalar equation \eqref{eq:1dconservationlaws} with the nonconvex flux $f(q)=\sin(q)$  and the initial condition
\begin{equation}
\label{sin_initial}
q_0(x) =\begin{cases}
\pi / 64, & \text{if}\     x<0,\\
255\pi/64, & \text{if}\     x\geq 0.
\end{cases}
\end{equation}
It is shown in Figure \ref{Fig:sin1}, that the numerical solution of the high order FV HWENO scheme does not converge to the entropy solution (solid black lines given by the first order Godunov scheme with a very refined mesh). One of the rarefaction waves in the compound wave is missing.


 \subsection{First order monotone modification}
\label{HWENO_mod1}

In this subsection, we propose a first order modification to the FV HWENO scheme for 1D nonconvex conservation laws following a similar idea in \cite{qiu_siam_sc}. The scheme can be summarized as follows, after a suitable initialization to obtain $\overline{q}^0$ and $\overline{\xi}^0$.

\begin{enumerate}
  \item Perform the HWENO reconstruction \cite{qiu2004hermite}.

  \noindent At each cell interface, say $x_{j+\frac12}$, reconstruct the point values $q_{j+\frac12}^\pm$ and derivative values $\xi_{j+\frac12}^\pm$ using neighboring cell average $\overline{q}$ and $\overline{\xi}$ respectively by the fifth order HWENO reconstruction procedure in Section \ref{HermiteWENO}.

  \item Identify the troubled cell boundary $x_{j+\frac12}$.

  \noindent\emph{ Criterion I: A cell boundary $x_{j+\frac12}$ is good, if $q_{j+\frac12}^\pm$, $\overline{q}_j$ and $\overline{q}_{j+1}$ all fall into the same linear, convex or concave region of the flux function $f(q)$. Otherwise, it is defined to be a troubled cell boundary. }

  \item At troubled cell boundaries, modify the numerical flux $\widehat{f}_{j+\frac12}$ and
  $\widehat{ \mathcal{H} }_{j+\frac12}$ with a discontinuity indicator in \cite{xu2006anti}.
      Specifically, the discontinuity indicator $\phi_j$ is defined as
      \begin{equation}
      \label{phi}
      \phi_j = \frac{  \beta_j}{\beta_j + \gamma_j}
      \end{equation}
  where
  \begin{equation*}
  \alpha_j=|\overline{q}_{j-1} - \overline{q}_j   |^2 + \varepsilon,\     \tau_j=|\overline{q}_{j+1}-\overline{q}_{j-1}|^2+\varepsilon,
 \      \beta_j = \frac{\tau_j}{\alpha_{j-1}} + \frac{\tau_j}{\alpha_{j+2}}, \      \gamma_j=\frac{ (q_{max}-q_{min})^2 }{\alpha_j}.
  \end{equation*}
  Here $\varepsilon$ is a small positive number taken as $10^{-6}$ in the code, and $q_{max}$ and $q_{min}$ are the maximum and minimum values of $\overline{q}_j$ over all cells. The discontinuity indicator $\phi_j$ has the property that

\begin{itemize}
  \item $0\leq \phi_j \leq1.$

  \item  $\phi_j$ is on the order of $O(\Delta x^2)$ in smooth regions.

  \item $\phi_j$ is close to $O(1)$ near a strong discontinuity.
\end{itemize}

Let $\widehat{f}_{j+\frac12} = \widehat{f}\big( q_{j+\frac12}^{m,-} , q_{j+\frac12}^{m,+}\big),$ and $\widehat{ \mathcal{H} }_{j+\frac12} = \widehat{  \mathcal{H} }\big( q_{j+\frac12}^{m,-} , q_{j+\frac12}^{m,+} ;\xi_{j+\frac12}^{m,-} ,\xi_{j+\frac12}^{m,+}\big)$ where
\begin{equation}
\label{eq: first_q}
q_{j+\frac12}^{m,-} = (1-\phi_j^2)q_{j+\frac12}^- + \phi_j^2\overline{q}_{j}, \quad q_{j+\frac12}^{m,+} = (1-\phi_j^2)q_{j+\frac12}^+ + \phi_j^2\overline{q}_{j+1},
\end{equation}
\begin{equation}
\label{eq: first_xi}
\xi_{j+\frac12}^{m,-} = (1-\phi_j^2)\xi_{j+\frac12}^- + \phi_j^2\overline{\xi}_{j}, \quad \xi_{j+\frac12}^{m,+} = (1-\phi_j^2)\xi_{j+\frac12}^+ +\phi_j^2\overline{\xi}_{j+1},
\end{equation}
with $\phi_j$ defined by \eqref{phi}, if $x_{j+\frac12}$ is a troubled cell boundary. Otherwise, at good cell boundaries, $q_{j+\frac12}^{m,\pm} = q_{j+\frac12}^{\pm}$ and $\xi_{j+\frac12}^{m,\pm} = \xi_{j+\frac12}^{\pm}$.

  \item  Evolve the cell averages $\overline{q}_j$ and $\overline{\xi}_j$ by \eqref{label:scheme}.

\end{enumerate}

\begin{remark}
When a troubled cell boundary is at a strong discontinuity, $\phi_j \sim 1$, hence $q_{j+\frac12}^{m,-} \sim \overline{q}_j$,   $q_{j+\frac12}^{m,+} \sim \overline{q}_{j+1},$ $\xi_{j+\frac12}^{m,-} \sim \overline{\xi}_j$ and $\xi_{j+\frac12}^{m,+} \sim \overline{\xi}_{j+1},$ indicating a first order monotone scheme is taking effect around a nonconvex discontinuous region. When a troubled cell boundary is in a smooth region, the modification is obtained with the magnitude at most of the size
\begin{equation*}
\phi_j^2 \max \big(\big|q_j-q_{j+\frac12}^-\big|,\big| q_{j+1}-q_{j+\frac12}^+\big|\big) \sim O(\Delta x^5),
\end{equation*}
\begin{equation*}
\phi_j^2 \max \big(\big|\xi_j-\xi_{j+\frac12}^-\big|,\big| \xi_{j+1}-\xi_{j+\frac12}^+\big|\big) \sim O(\Delta x^5),
\end{equation*}
hence it does not affect the fifth order accuracy of the scheme.
\end{remark}

It is natural to extend the above first order modification to  2D problems. With the system, the HWENO reconstructions are performed in local characteristics directions \cite{qiu2004hermite}, then a first order monotone modification in the form of \eqref{eq: first_q} and \eqref{eq: first_xi} is applied. For 2D problems, we identify trouble cell boundaries, that is to check if the convexity fails, via Gaussian points along cell boundaries, e.g. $q_{i+\frac12,i_g}^\pm$ in equation \eqref{eq: flux_r}. Similarly the first order modification is performed with respect to these Gaussian points along cell boundaries.

\subsection{Second order modification with an entropic projection}
\label{second}
A MUSCL type method with an entropic projection is proposed in \cite{bouchut1996muscl}. It is proved in the same paper that schemes with such entropic projection enjoy cell entropy inequalities for all convex entropy functions. In the following, we apply such entropic projection as a modification to the HWENO scheme around nonconvex regions to ensure entropic convergence.

\subsubsection{Review of the MUSCL method satisfying all the numerical entropy inequalities.}
The procedure of the MUSCL scheme satisfying all entropy inequalities can be summarized below.
Let the numerical solution at time level $n$ be written as $q^n_j(x) = \overline{q}_j^n + s_j^n \sigma_j$ with $\sigma_j = \frac{x-x_j}{\Delta x_j}$ over the cell $I_j$.
Initially, $q^0_j(x) = \overline{q}_j^0 + s_j^0 \sigma_j$, with
$
s_j^0 = \text{minmod}( \overline{q}_j^0-\overline{q}_{j-1}^0, \overline{q}_{j+1}^0-\overline{q}_j^0 )
$
where the minmod function is defined as follows,
\begin{equation*}
\text{minmod}(a,b)=
\begin{cases}
0, & \text{if}\ ab<0,
\\
 \min(a,b), &
\ \text{if}\  a,b\geq 0,
\\
 \max(a,b), &
\ \text{if} \ a,b\leq 0.
\end{cases}
\end{equation*}
It consists of two steps to evolve from $q^n$ to $q^{n+1}$.

\begin{enumerate}
\item {\emph{Exact evolution}} $(T_{\Delta x})$: Evolve \eqref{eq:1dconservationlaws} exactly for a time step $\Delta t$, to obtain a solution $\widetilde{q}^{n+1}$, which in general is not a piecewise linear function anymore.

    \item \emph{An entropic projection} $(P^1)$: Find a second order approximation to $\widetilde{q}^{n+1}$ by a piecewise linear function $q^{n+1}$, satisfying
        \begin{equation}
        \label{entropy}
        \int_{I_j}U(q^{n+1}(x))dx \leq \int_{I_j} U (\widetilde{q}^{n+1}(x)) dx, \    \forall j
        \end{equation}
        for all convex entropy function $U(u)$. Second order reconstruction satisfying \eqref{entropy} can be obtain by setting the cell average as
        \begin{equation}
        \overline{q}_j^{n+1} = \frac{1}{\Delta x_j}\int_{I_j} \widetilde{q}^{n+1}(x) dx
        \end{equation}
        and the slope as
        \begin{equation}
        s_j^{n+1} = D \widetilde{q}^{n+1}|_{I_j} = \text{minmod}_{I_j} \zeta(y)
        \end{equation}
        where
        \begin{equation}
        \label{zeta}
        \zeta(y) = \frac{2}{\Delta x_j}\left(  \frac{1}{x_{j+\frac12}-y}\int_y^{x_{j+\frac12}}\widetilde{q}^{n+1}(x) dx - \frac{1}{y-x_{j-\frac12}}\int_{x_j-\frac12}^y \widetilde{q}^{n+1}(x)dx  \right)
        \end{equation}
        The minmod function of $g(x)$ on the interval $(a,b)$ is defined as
        \begin{equation}
        \mathop{ \text{minmod} }_{(a,b)} g(x) =
        \left\{
        \begin{array}{ll}
        0,   &   \text{if} \    \exists y_1,y_2 \in(a,b),\    s.t.   \   g(y_1)g(y_2)\leq 0, \\
        \mathop{\min}_{(a,b)} g(y),   &   \text{if} \    g(y) >0,\     \forall y\in(a,b),\\
        \mathop{\max}_{(a,b)} g(y),   &   \text{if} \    g(y) <0, \    \forall y\in(a,b).
        \end{array}
        \right.
        \end{equation}
\end{enumerate}

In summary, the scheme can be written out in the following abstract form
\begin{equation}
\label{bouchut1996muscl}
q^{n+1}=P^1 \circ T_{\Delta t}(q^n) \doteq Q^1(\Delta t)(q^n).
\end{equation}

It enjoys the following convergence theorem as proved in \cite{bouchut1996muscl}.

\begin{theorem}
\cite{bouchut1996muscl}. Let $T=n\Delta t$, $u(\cdot,T)$ be the exact entropy solution to \eqref{eq:1dconservationlaws} with the initial data $q_0\in L^1 \cap BV(\mathcal{R}),f_\infty = \max_{q\in [\min q_0,\max q_0]}f'(q),$ then there exists a constant $C$, such that
\begin{equation}
\| Q^1(\Delta t)^n q^0 - q(\cdot,T) \|_{L^1} \leq C(f_\infty \sqrt{T\Delta t} + \Delta x\sqrt{T/\Delta t}) .
\end{equation}
The second order MUSCL scheme with the entropic projection \eqref{bouchut1996muscl} converges to the unique entropy solution, when $\Delta t = \mathcal{O}(\Delta x)$.
\end{theorem}

\subsubsection{Second order modification to the fifth order FV HWENO schemes}

At each time step evolution, $\{  \overline{q}_j^\star , \overline{\xi}_j^\star, q_j^{\star,l},q_j^{\star,r} \}$, $\star = n,(1),(2)$, over the cell $I_j$ is updated in each RK stage. For instance, at the initial stage, $q_{j\mp \frac12}^{n,\pm}$ is obtain by the HWENO reconstruction from $\overline{q}^{n}$ and $\overline{\xi}^{n}$. $q_j^{n,l}$ and $q_j^{n,r}$ refer to approximations with {\em entropic projection} to the left and right boundaries of $I_j$ when the cell is detected as a trouble cell. Initially, $q_j^{n,l}$ and $q_j^{n,r}$ come from a MUSCL scheme with a minmod reconstruction.
To show the idea of second order modification to the fifth order FV HWENO schemes, we present the procedure to update $\{  \overline{q}_j^{(1)} , \overline{\xi}_j^{(1)},   q_j^{(1),l},q_j^{(1),r} \}$ from  $\{  \overline{q}_j^n , \overline{\xi}_j^n,   q_j^{n,l},q_j^{n,r} \}$.

\begin{description}
  \item[Step 1.] Compute $q_{j\pm\frac12}^{n,\pm} $ by performing the HWENO reconstruction from $\{ \overline{q}_j^n , \overline{\xi}_j^n \}$.

  \item[Step 2.] Update $\overline{q}_j^{(1)}$ and $\overline{\xi}_j^{(1)}$:
  \begin{enumerate}
  \item Identify the troubled cell boundaries, for which we refer to Criterion I in section \ref{HWENO_mod1} for the details.

  \item Only at trouble cell boundaries, modify numerical fluxes $\widehat{f}_{j+\frac12}$ and $\widehat{ \mathcal{H}  }_{j+\frac12}$ as follows. We let $\widehat{f}_{j+\frac12} = \widehat{f}\big(u_{j+\frac12}^{m,-},u_{j+\frac12}^{m,+}\big)$ and $\widehat{ \mathcal{H} }_{j+\frac12} = \widehat{ \mathcal{H} }\big(u_{j+\frac12}^{m,-},u_{j+\frac12}^{m,+} ; \xi_{j+\frac12}^{m,-},\xi_{j+\frac12}^{m,+}   \big)$, where
  \begin{equation}
  q_{j+\frac12}^{m,-} = (1-\phi_j^2)q_{j+\frac12}^{n,-} + \phi_j^2 q_j^{n,r}, \quad q_{j+\frac12}^{m,+}=(1-\phi_j^2)q_{j+\frac12}^{n,+} +\phi_j^2 q_{j+1}^{n,l},
  \end{equation}
    \begin{equation}
  \xi_{j+\frac12}^{m,-} = (1-\phi_j^2)\xi_{j+\frac12}^{n,-} + \phi_j^2 \overline{\xi}_j^n, \quad \xi_{j+\frac12}^{m,+}=(1-\phi_j^2)\xi_{j+\frac12}^{n,+} +\phi_j^2 \overline{\xi}_{j+1}^n,
  \end{equation}
  with $\phi_j$ defined by \eqref{phi} at the troubled cell boundary.

  \item Update the solution at the first RK stage as follows,
  \begin{equation*}
  \overline{q}_j^{(1)} = \overline{q}_j^n - \frac{\Delta t}{\Delta x} ( \widehat{f}_{j+\frac12} - \widehat{f}_{j-\frac12} ),
  \end{equation*}
   \begin{equation*}
  \overline{\xi}_j^{(1)} = \overline{\xi}_j^n - \frac{\Delta t}{\Delta x} ( \widehat{ \mathcal{H} }_{j+\frac12} - \widehat{ \mathcal{H} }_{j-\frac12} ).
  \end{equation*}
  \end{enumerate}

  \item[Step 3.]  Update $q_j^{(1),l}$ and $q_j^{(1),r}$ for a nonconvex troubled cell $I_j$: we first identify nonconvex troubled cells by the following criterion.

  \emph{Criterion II: A cell $I_j$ is called a good cell, if $q^m=\{ \overline{q}_m^n,q_{m\mp\frac12}^{n,\pm},q_m^{n,l},q_m^{n,r} \}$ with $m=j-1,j,j+1$, fall into the same linear, convex or concave region of the flux function $f(q)$. Otherwise, it is defined to be a nonconvex troubled cell. }

  \begin{description}
    \item[Trouble cells.]
  At a nonconvex troubled cell $I_j$,
  we apply a first order scheme on a refined mesh by evolving a time step $\Delta t$.
    Specifically, we evolve equation \eqref{eq:1dconservationlaws} with the initial condition
  \begin{equation}
  \overline{q}_l^n + s_l^n \sigma_l,\   \text{for} \   x\in I_l, \   l=j-1,j,j+1
  \end{equation}
  where $s_l^n = 2 minmod (q_l^{n,r} - \overline{q}_l^n, \overline{q}_l^n-q_l^{n,l})$ and $\sigma_l = \frac{x-x_l}{\Delta x_l}$.
  A periodic boundary condition on $I_{j-1}\cup I_j \cup I_{j+1}$ is used. We consider a refined numerical mesh of cell $I_j$
  \begin{equation}
  I_j = \cup_{m=1}^N [y_{m-\frac12},y_{m+\frac12}],\    \delta x = y_{m+\frac12} - y_{m-\frac12} = \Delta x/N,
  \end{equation}
  and apply a first order scheme with entropic convergence to evolve the solution for $\Delta t$. Let $\widetilde{q}^{(1)}|_{I_j}$ be the evolved solution, approximated by a piecewise constant function sitting on the refined numerical mesh with the truncation error $\sim \mathcal{O}(\delta x)=\mathcal{O}(\Delta x^2)$. We compute the average and slope of the linear function approximating $\widetilde{q}^{(1)}|_{I_j}$ on $I_j$ via the entropic projection as follows: the average is taken as the average of $\widetilde{q}^{(1)}|_{I_j}$ and the slope is computed as follows
        \begin{equation}
        s_j^{(1)} = \text{minmod} ( \zeta(y_{N+\frac12} ), \cdots, \zeta(y_{2N+\frac12} ) ) ,
        \label{expensive}
        \end{equation}
        where
        \begin{equation}
        \zeta(y) =
        \frac{2}{\Delta x_j}
        \left(
         \frac{1}{ x_{j+\frac12} -y } \int_{ y }^{x_{j+\frac12} }  \widetilde{q}^{(1)} dx -
         \frac{1}{ y -x_{j-\frac12} } \int_{ x_{j-\frac12} }^y \widetilde{q}^{(1)} dx
        \right)
        \end{equation}

        Finally,
        \begin{equation}
        q_j^{ (1),l } = \overline{ \widetilde{q}_j^{(1)}  }  - \frac12 s_j^{(1)}, \
        q_j^{ (1),r } = \overline{ \widetilde{q}_j^{(1)}  }  + \frac12 s_j^{(1)}.
        \end{equation}

    \item[Good cells.]
 For a good cell $I_j$, update $q_j^{(1),l}$ and $q_j^{(1),r}$ by setting $q_j^l = q_{j-\frac12}^{(1),+}$ and $q_j^r = q_{j+\frac12}^{(1),-}$, where $q_{j-\frac12}^{(1),+}$ and $q_{j+\frac12}^{(1),-}$ are reconstructed by performing a HWENO reconstruction.
\end{description}

 \end{description}

\begin{remark}
Note that the modification of $\xi_{j+\frac12}^\pm$ is a first order modification on derivative values. Because the first order on derivative values is enough to get a second order scheme.
\end{remark}

\begin{remark}
The implementation of the procedure to find $s_j^{n+1}$ in \eqref{expensive} is computationally expensive.
Due to the costly implementation of the high order scheme with the second order entropic projection, we only adopt the first order modification to modify the FV HWENO for   2D scalar problems.
\end{remark}

\section{Numerical Experiments}
\label{numericaltests}

In Section \ref{1D_scalar}, we compare the performance of the fifth order FV HWENO scheme (HWENO5) and the fifth order FV WENO scheme (WENO5) with the first order modification (mod1) and the second order entropic projection (mod2) respectively for solving 1D nonconvex conservational laws.
In Section \ref{2D_scalar}, we present the performance of the modified WENO5 and HWENO scheme (HWENO4) for 2D problems. The numerical fluxes used in this paper are the global Lax-Friedrich flux.

\subsection{1D scalar problems}
\label{1D_scalar}

\begin{example}
We consider the nonconvex conservation law
\begin{equation}
q_t + \left( \frac{q^3}{3}\right)_x = 0, \     q_0(x) = \sin(\pi x).
\label{u3}
\end{equation}
We compute the solution up to $t=0.2$. Table \ref{table:11} gives the $L^1$ and $L^\infty$ errors and the corresponding orders of accuracy of the regular and modified HWENO5 and WENO5 schemes. We can  see that errors of HWENO5
are smaller than those of WENO5 with the same mesh. Very little difference is observed among the regular and two modified HWENO5 and WENO5 schemes.

\begin{table}[!ht]\small
\caption{\small $q_t + \left( \frac{q^3}{3}\right)_x = 0$ with initial condition $q_0(x) = \sin(\pi x)$ and periodic boundary conditions. The $L^1$ and $L^\infty$ errors and the corresponding orders of accuracy for the regular HWENO5 and WENO5, the corresponding two versions of modified schemes at the time $t=0.2$. }
\vspace{0.1in}
\centering
\begin{tabular}{|l||llll||llll|}
\hline
{\multirow{2}{*}{N}} & \multicolumn{4}{c||}{HWENO5}  &\multicolumn{4}{c|}{WENO5}   \\         
\multicolumn{1}{|l||}{}&
 \multicolumn{1}{|l}{$L_1$ error} & Order & $L_\infty$ error & Order & {$L_1$ error} & Order & $L_\infty$ error & Order \\
\hline
 & \multicolumn{8}{c|}{ regular}  \\
\hline
   100 &     1.68E-05 & &     2.14E-04 &     &     4.42E-05 & &     4.61E-04 &  \\ \hline
   200 &     7.68E-07 &     4.45 &     1.69E-05 &     3.66  &     2.24E-06 &     4.30 &     4.59E-05 &     3.33 \\ \hline
   300 &     1.17E-07 &     4.64 &     2.62E-06 &     4.59 &     3.48E-07 &     4.59 &     7.68E-06 &     4.41  \\ \hline
   400 &     2.96E-08 &     4.78 &     6.44E-07 &     4.88 &     9.00E-08 &     4.71 &     1.96E-06 &     4.76  \\ \hline
   500 &     1.01E-08 &     4.81 &     2.12E-07 &     4.97 &     3.12E-08 &     4.75 &     6.58E-07 &     4.88  \\ \hline
    & \multicolumn{8}{c|}{mod1}  \\
\hline
   100 &     1.68E-05 & &     2.14E-04 &     &     4.42E-05 & &     4.61E-04 &  \\ \hline
   200 &     7.68E-07 &     4.45 &     1.69E-05 &     3.66  &     2.24E-06 &     4.30 &     4.59E-05 &     3.33 \\ \hline
   300 &     1.17E-07 &     4.64 &     2.62E-06 &     4.59 &     3.48E-07 &     4.59 &     7.68E-06 &     4.41  \\ \hline
   400 &     2.96E-08 &     4.78 &     6.44E-07 &     4.88 &     9.00E-08 &     4.71 &     1.96E-06 &     4.76  \\ \hline
   500 &     1.01E-08 &     4.81 &     2.12E-07 &     4.97 &     3.12E-08 &     4.75 &     6.58E-07 &     4.88  \\ \hline
    & \multicolumn{8}{c|}{mod2}  \\
\hline
   100 &     1.70E-05 & &     2.14E-04 &   &        4.42E-05 & &     4.61E-04 & \\ \hline
   200 &     7.78E-07 &     4.45 &     1.69E-05 &     3.66   &     2.24E-06 &     4.30 &     4.59E-05 &     3.33 \\ \hline
   300 &     1.17E-07 &     4.67 &     2.62E-06 &     4.59   &     3.48E-07 &     4.59 &     7.68E-06 &     4.41 \\ \hline
   400 &     2.96E-08 &     4.78 &     6.44E-07 &     4.88  &     9.00E-08 &     4.70 &     1.96E-06 &     4.76 \\ \hline
   500 &     1.01E-08 &     4.81 &     2.12E-07 &     4.97   &     3.12E-08 &     4.75 &     6.58E-07 &     4.88 \\ \hline
\end{tabular}
\label{table:11}
\end{table}

\end{example}

\begin{example}
Consider the Riemann problem of the nonconvex conservation law presented in Section~\ref{ssub: bad_example}.
We plot numerical solutions of two modified schemes in Fig.~\ref{Fig:sin_mod}. They both successfully converge to the correct entropy solution, with the development of a compound wave including a shock, a rarefaction wave, followed by another shock and another rarefaction wave.
\begin{figure}[h!]
\centering                              
\includegraphics[height=60mm]{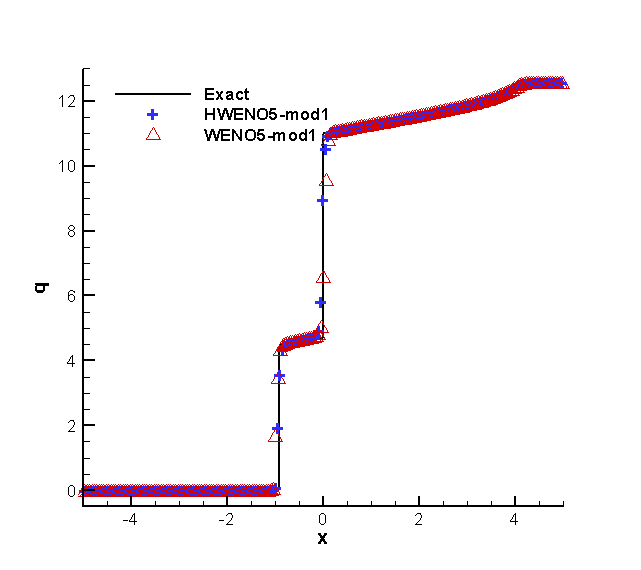}
\includegraphics[height=60mm]{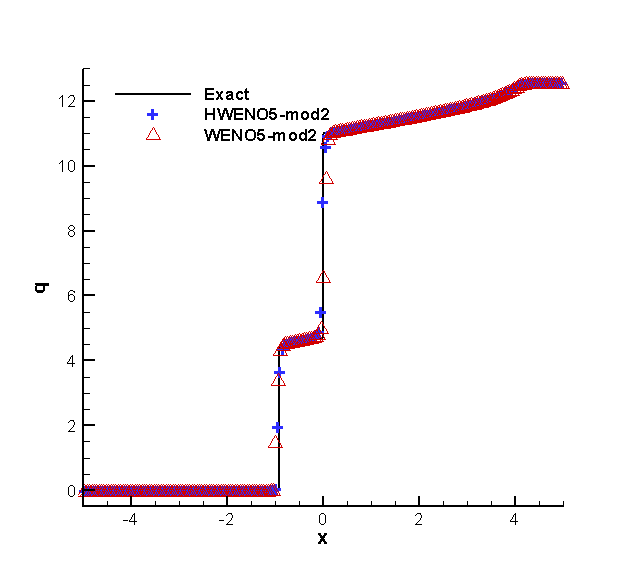}
\caption{\small Solid lines: the exact solution of \eqref{eq:1dconservationlaws} with $f(u)=\sin(u)$ and the initial condition \eqref{sin_initial} at the time $t=4$; HWENO5 (pluses) and  WENO5 (deltas) with the first order monotone schemes (left); HWENO5 (pluses) and WENO5 scheme (deltas) with the second order entropy projection (right). using $N = 200$ uniform cells. }
\label{Fig:sin_mod}
\end{figure}
\end{example}

\begin{example}
Consider  \eqref{eq:1dconservationlaws}
 with the nonconvex flux $f(q)$ defined by
\begin{equation}
\label{nonconvex}
f(q)=\left\{
\begin{array}{ll}
1, & \text{if} \quad  q<1.6 \\
\cos(5\pi(q-1.8)) +2.0, & \text{if} \quad 1.6 \leq q <2.0 \\
-\cos(5\pi(q-2.2)),   & \text{if} \quad 2.0\leq q <2.4   \\
1,                    &    \text{if}  \quad q\geq 2.4
\end{array}
 \right.
\end{equation}
with the initial condition
\begin{equation}
\label{init1}
q_0(x) = \left\{
\begin{array}{ll}
1, & \text{for} \quad x<0\\
3, & \text{for} \quad x\geq0.
\end{array}\right.
\end{equation}
In the left panel of Figure \ref{Fig:eg1mod}, the HWENO5 seems to converge to the entropy solution slowly, which might be related to the fact that the reconstruction of the solution at the rarefaction waves comes from neighboring cells and is not a good approximation when the rarefaction wave is surrounded by two shocks at its early stage of development.
As shown in Figure  \ref{Fig:eg1mod}, the numerical solutions of two modified schemes successfully converge to the correct entropy solution.

\begin{figure}[h!]
\centering                              
\includegraphics[height=60mm]{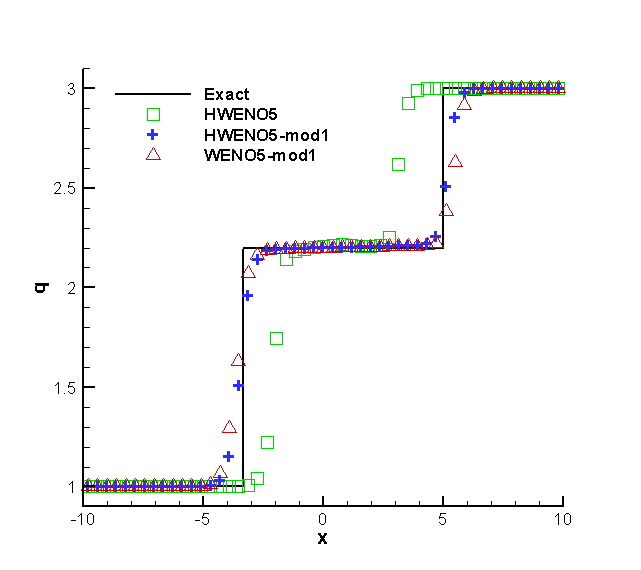}
\includegraphics[height=60mm]{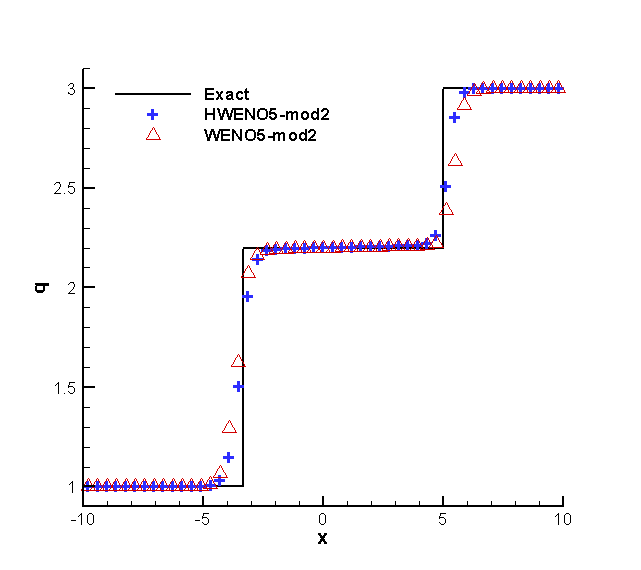}
\caption{\small Solid line: The exact solutions of the nonconvex scalar conservation law \eqref{eq:1dconservationlaws}-\eqref{nonconvex} with the initial condition \eqref{init1} at the time $t=2$. Left: HWENO5 (squares), HWENO5-mod1 (pluses) and WENO5-mod1 (deltas); Right: HWENO5-mod2 (pluses) and WENO5-mod2 (deltas).
$N = 50$ uniform cells are used. CFL = 0.5.}
\label{Fig:eg1mod}
\end{figure}

\end{example}

\begin{example}
The nonconvex conservation law \eqref{eq:1dconservationlaws} with the nonconvex flux $f(q)$ given by \eqref{nonconvex} with the initial condition
\begin{equation}
q_0(x) =
\begin{cases}
3, & \text{for} \     -1\leqslant x<0 \\
1, &  \text{for} \     0\leqslant x\leqslant 1
\end{cases}
\label{31}
\end{equation}
and a periodic boundary condition. This is a very challenging test case: with periodic boundary conditions, the compound waves strongly interact with each other.
There is no analytic formula of the exact solution for this problem. The reference solution is computed by the Godunov scheme with 400,000 uniform cells.
It is observed from Fig.~\ref{Fig:31_WL} that, numerical solutions of the HWENO5 scheme without modification deviate away from the reference solution with mesh refinement.
For this example, due to strong interaction of compound waves, the monotonicity preserving limiter (MPHWENO5) \cite{suresh1997accurate} is applied to control oscillations.
As shown in Fig \ref{Fig:31_converge}, numerical solutions of modified HWENO5 schemes converge to the correct entropy solution.
The comparison of WENO5/MPHWENO5 with different modifications is shown in Fig \ref{Fig:31_LF}.
%
The numerical solution of HWENO5 with the first order modification is observed to converge to the correct entropy solution slightly faster, compared to that of WENO5 with the first order modification.
Comparable performance of  HWENO5 and WENO5 with the second order modification are observed.
We observe better performance of HWENO5 or WENO5 with the second order modification when compared with schemes with a first order modification.
\begin{figure}[h!]
\centering                              
\includegraphics[height=70mm]{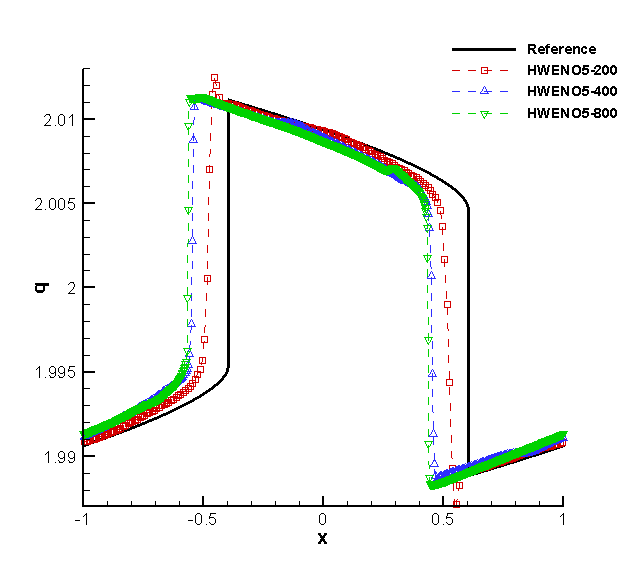}
\caption{\small Solid lines: The reference solution of \eqref{31} at the time $t=2$;  HWENO5 without modification with $N=200$ uniform cells (Squares), with $N=400$ uniform cells (Deltas) and with $N=800$ uniform cells (Gradients). CFL = 0.01.}
\label{Fig:31_WL}
\end{figure}

\begin{figure}[h!]
\centering                              
\includegraphics[height=70mm]{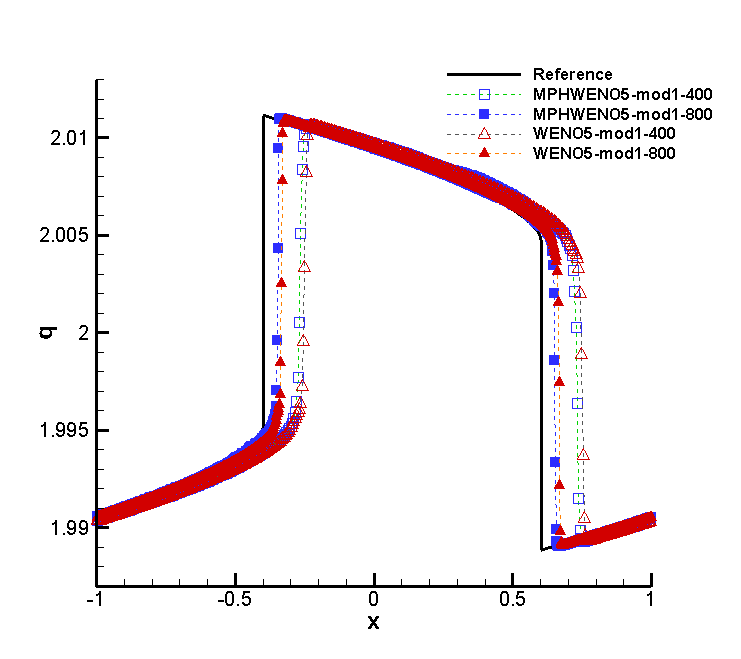}
\caption{\small Solid lines: The reference solution of \eqref{31} at the time $t=2$;  MPHWENO5 with the first order modification with $N=400$ uniform cells (Squares) and with $N=800$ uniform cells (Filled squares);  WENO5 with the first order modification with $N=400$ uniform cells (Deltas) and with $N=800$ uniform cells (Filled deltas). CFL = 0.01.}
\label{Fig:31_converge}
\end{figure}

\begin{figure}[h!]
\centering                              
\includegraphics[height=70mm]{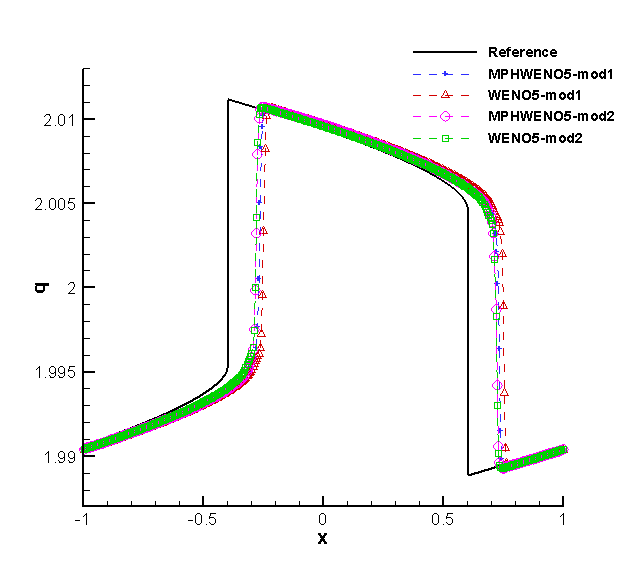}
\includegraphics[height=70mm]{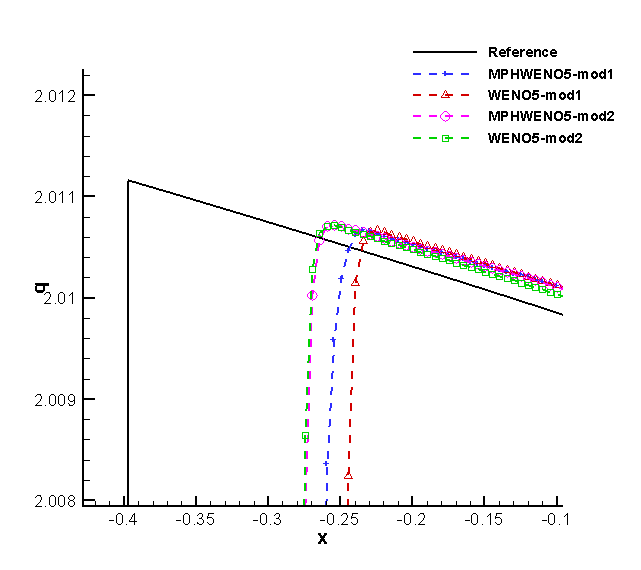}
\caption{\small Solid lines: The reference solution of \eqref{31} at the time $t=2$; Pluses: MPHWENO5 with the first order modification; Deltas: WENO5 with the first order modification; Circles: MPHWENO5 with the second order modification; Squares: WENO5 with the second order modification; The zoom are given in the right; $N = 400$ uniform cells are used; CFL = 0.01.}
\label{Fig:31_LF}
\end{figure}

\end{example}

\subsection{2D scalar problems}
\label{2D_scalar}
\begin{example}
We solve the following nonconvex conservation law in 2D :
\begin{equation*}
q_t + \left(\frac{q^3}{3}\right)_x + \left(\frac{q^3}{3}\right)_y =0,
\end{equation*}
with the initial condition $q(x,y,0)=\sin(\pi(x+y)/2)$ and the periodic boundary condition
in both directions. The computational
domain for this problem is $[-2,2]\times[-2,2]$. When $t=0.2$ the solution is still smooth.


\begin{table}[!ht]\small
\caption{\small
$q_t + \left(\frac{q^3}{3}\right)_x + \left(\frac{q^3}{3}\right)_y =0$ with initial condition $q(x,y,0)=sin(\pi(x+y)/2)$ and periodic boundary conditions.
 The $L^1$ and $L^\infty$ errors and the corresponding orders of accuracy for the regular HWENO4 and WENO5, the corresponding modified schemes at the time $t=1$.}
\vspace{0.1in}
\centering
\begin{tabular}{|l||llll||llll|}
\hline
{\multirow{2}{*}{N}} & \multicolumn{4}{c||}{HWENO4}  &\multicolumn{4}{c|}{WENO5}   \\         
\multicolumn{1}{|l||}{}&
 \multicolumn{1}{|l}{$L_1$ error} & Order & $L_\infty$ error & Order & {$L_1$ error} & Order & $L_\infty$ error & Order \\
\hline
 & \multicolumn{8}{c|}{ regular}  \\
\hline
40  &    3.56E-03 &    2.36 &       6.83E-03 &     0.72 &    5.26E-04 &     2.63  &       4.23E-03 &     1.08
    \\ \hline
80 &    4.25E-04 &    3.07 &       1.20E-03 &     2.51 &    7.73E-05 &     2.77  &       7.86E-04 &     2.43
   \\ \hline
160 &    2.72E-05 &    3.96 &       1.15E-04 &     3.39 &    5.35E-06 &     3.85 &       8.29E-05 &     3.25
   \\ \hline
320 &    1.16E-06 &    4.56 &       6.65E-06 &     4.11 &   2.31E-07 &     4.54 &       4.58E-06 &     4.18
    \\ \hline
640 & 5.68E-08 &     4.35  &       3.16E-07 &     4.39  &   9.68E-09 &     4.57  &      1.39E-07 &     5.05
\\  \hline
    & \multicolumn{8}{c|}{ mod1}  \\
\hline
40   &    4.26E-03 &    2.51 &        6.84E-03 &     1.56&    5.95E-04 &     3.04 &        4.28E-03 &     1.97
    \\ \hline
80 &    4.82E-04 &    3.15 &        1.20E-03 &     2.51 &    7.97E-05 &     2.90 &        7.86E-04 &     2.44
   \\ \hline
160  &    3.11E-05 &    3.95 &        1.15E-04 &     3.39&    5.42E-06 &     3.88 &        8.29E-05 &     3.25
   \\ \hline
320  &    1.38E-06 &    4.50 &        6.65E-06 &     4.11&    2.32E-07 &     4.54 &        4.58E-06 &     4.18
    \\ \hline
640  &   6.99E-08 &     4.30  &       3.16E-07 &     4.39 &  9.73E-09 &     4.58  &        1.39E-07 &     5.05
 \\ \hline
\end{tabular}
\label{table:23}
\end{table}

Table \ref{table:23} gives the $L_1$ errors and the $L_{\infty}$ errors and the corresponding orders of the accuracy of the regular and modified FV HWENO scheme and FV WENO scheme. $\Delta t = \frac{CFL}{\frac{\alpha}{\Delta x} + \frac{\beta}{\Delta y} }$ where $CFL=0.4$, $\alpha=\mathop{\max}\{f'(q)\},\beta=\mathop{\max}\{g'(q)\}$.  Expected orders of convergence are observed.

\end{example}


\begin{example}
We solve the KPP rotating
wave problem,
\begin{equation*}
u_t + (\sin(u) )_x + (\cos(u) )_y =0
\end{equation*}
with the initial condition
\begin{equation*}
u(x,y,0)=
\left\{
\begin{array}{ll}
\frac{14}{4}\pi, &  \text{if} \   \sqrt{x^2+y^2}\leq 1, \\
\frac{\pi}{4},  &   \text{otherwise}.
\end{array}
\right.
\end{equation*}
This test was originally proposed in Kurganov et al.
\cite{kurganov2007adaptive}. It is challenging to many high-order numerical schemes because the solution has
a two-dimensional composite wave structure.

\begin{figure}[h!]
\centering                              
\includegraphics[height=65mm]{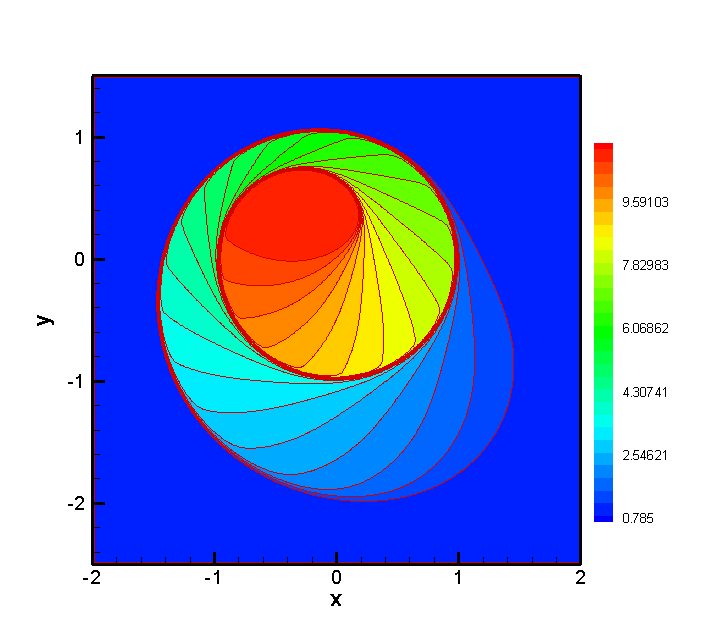}
\includegraphics[height=65mm]{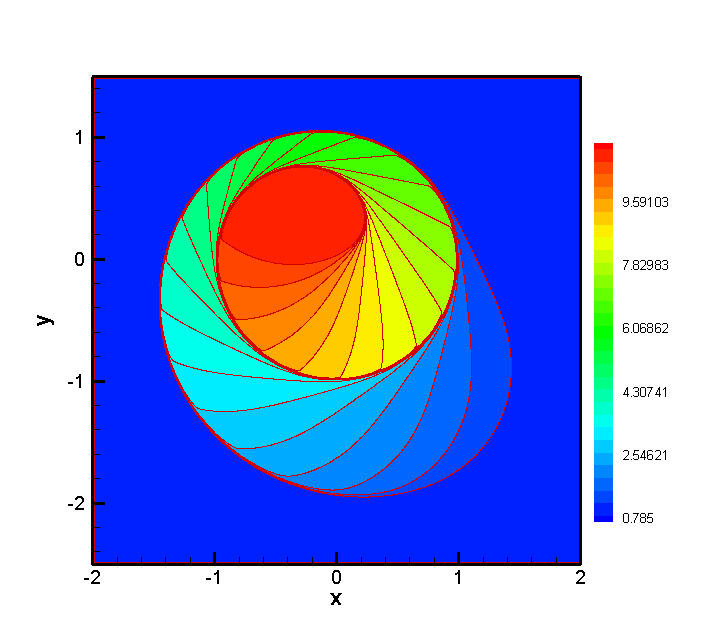}
\includegraphics[height=65mm]{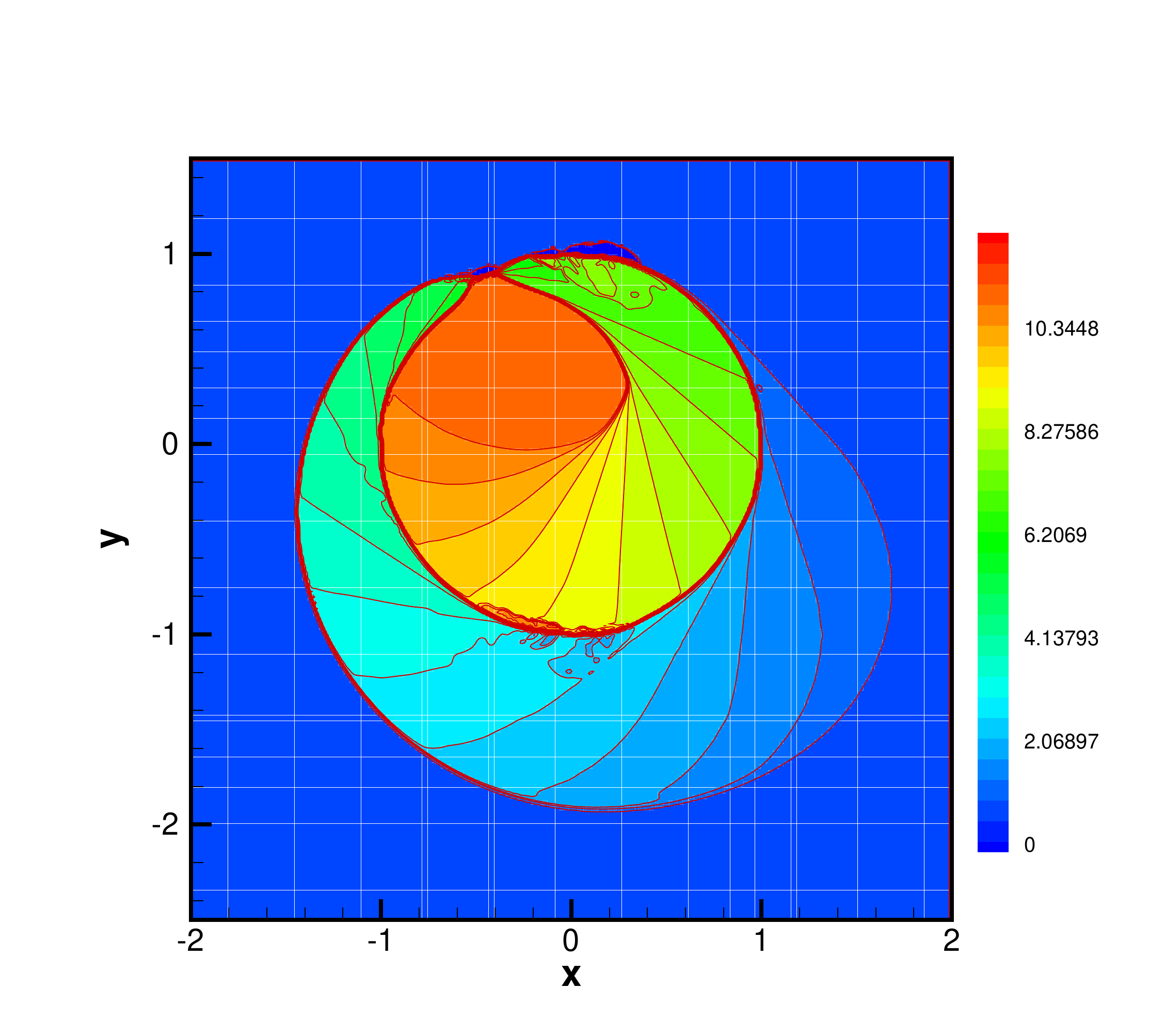}
\includegraphics[height=65mm]{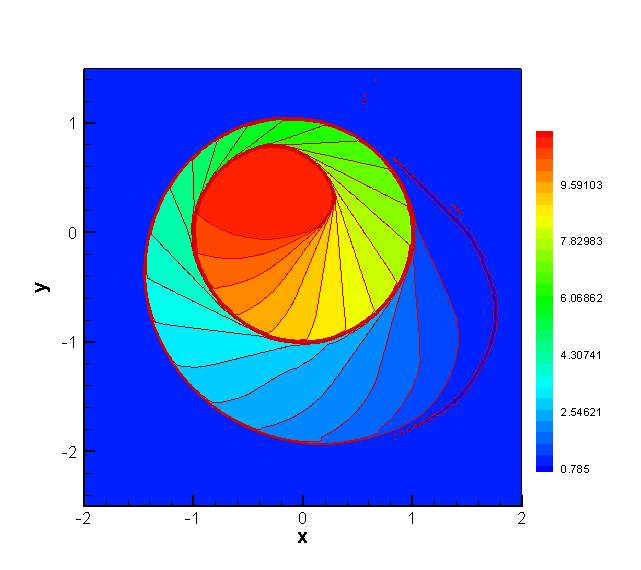}
\includegraphics[height=65mm]{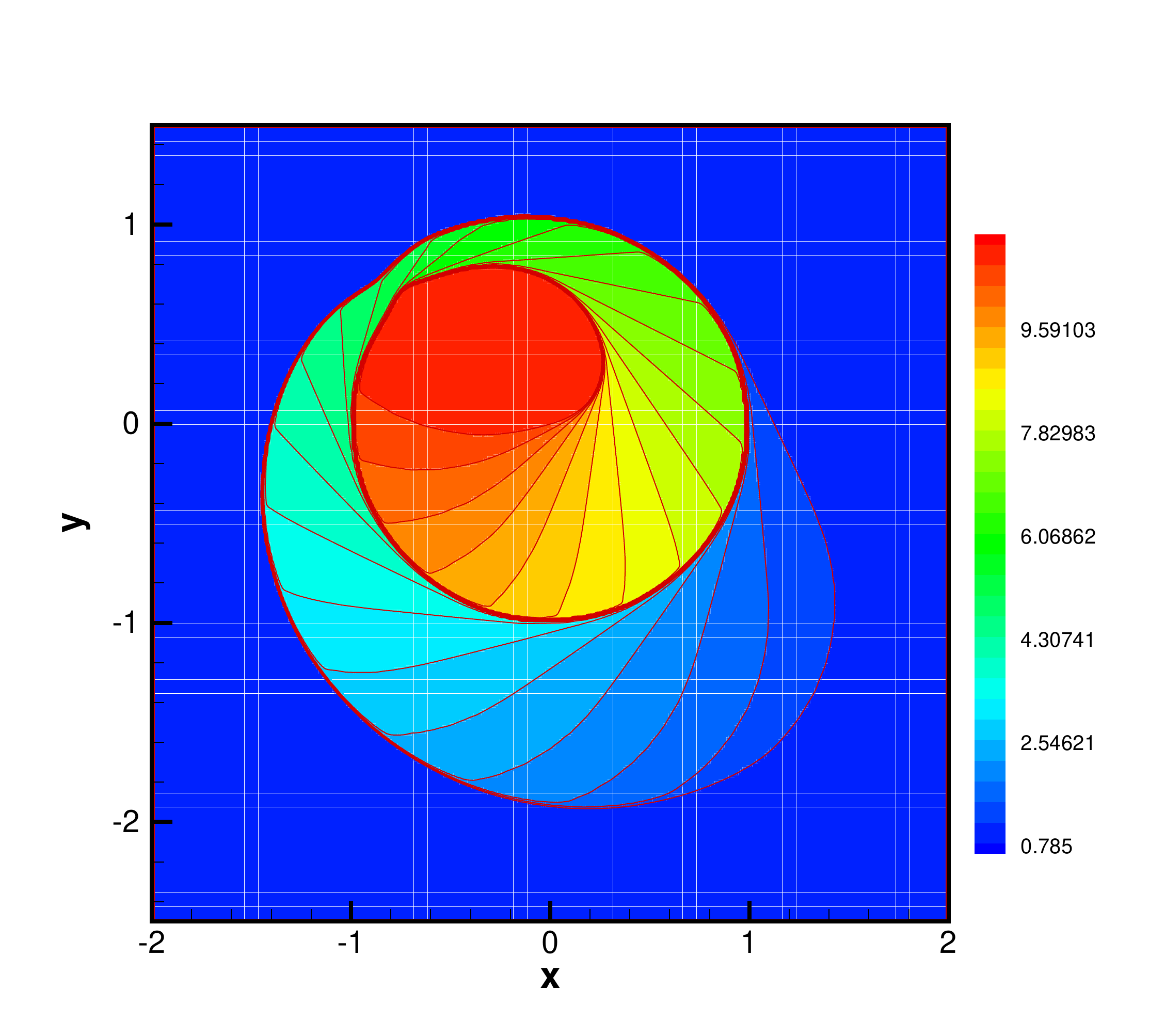}
\includegraphics[height=65mm]{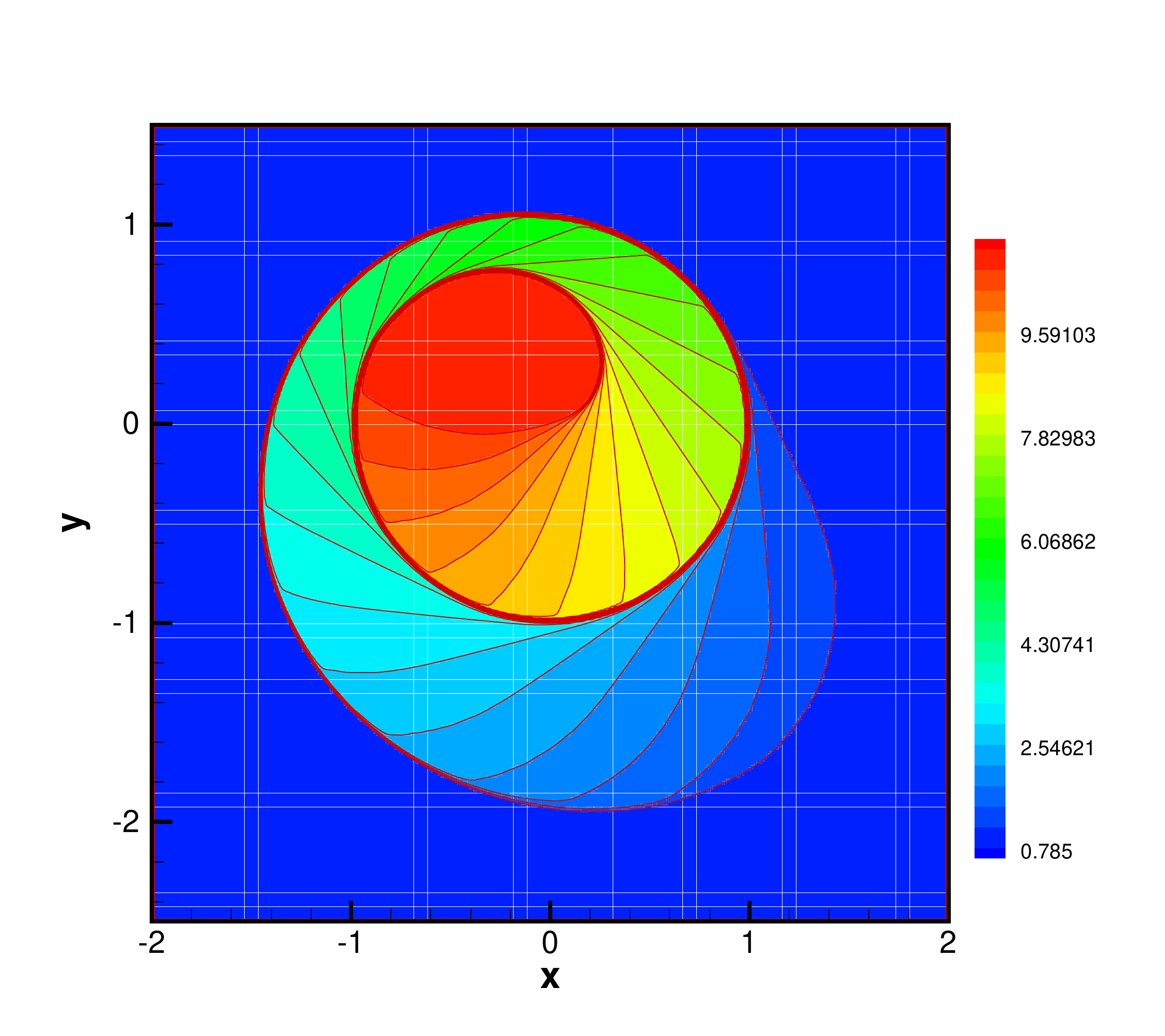}
\caption{\small KPP problem at time $t=1$. 30 equally
spaced solution contours from 0.785 to 11.0.
Fist row: first-order approximation with $400\times400$ cells; first-order approximation with $1000\times1000$ cells.
Second row: regular HWENO4 with $400\times400$ cells; modified HWENO4 with $400\times400$ cells.
Third row: regular WENO with $400\times400$ cells; modified WENO with $400\times400$ cells. }
\label{Fig:eg1}
\end{figure}

In Figure \ref{Fig:eg1}, we show the contours of the solution at $t=1$. In the left panel, it is observed that neither HWENO4 or WENO5 schemes can capture composite wave structures. The composite wave structures are well captured by the HWENO4 or WENO5 with the first order modification as shown on the right panel.

\end{example}

\section{Concluding remarks}
\label{concluding}
We proposed modifications to FV HWENO schemes for nonconvex conservation laws based on the idea of \cite{qiu_siam_sc}, emphasizing convergence to the entropy solution. The robustness of modified FV HWENO schemes is showed by several representative examples including  2D problems. We also compare the performance between the modified FV HWENO and WENO schemes.

\bibliographystyle{siam}
\bibliography{refer}

\end{document}